\documentclass[a4paper,11pt]{article}
\usepackage{graphicx,multirow}
\usepackage{amssymb,amsmath,amsthm}
\usepackage{lscape,lipsum,microtype}
\usepackage[margin=2.54cm]{geometry}  
\usepackage{hyperref}

\usepackage{ifluatex}
\ifluatex
 \usepackage{unicode-math}
\fi

\usepackage{latexsym}
\usepackage{amsmath}
\usepackage{amsfonts}
\usepackage{amssymb}
\usepackage{amscd}
\usepackage{amsmath}

\newcommand{\D }{\Delta }

\renewcommand{\l }{\lambda }

\newcommand{\n }{\nabla }

\newcommand{\intbar}{\mathop{\int\makebox(-13.5,0){\rule[4pt]{.7em}{0.3pt}}%
\kern-6pt}\nolimits}

\newcommand{\be}{\begin{equation}}
\newcommand{\ee}{\end{equation}}
\newcommand{\bes}{\begin{equation*}}
\newcommand{\ees}{\end{equation*}}
\newcommand{\ba}{\begin{eqnarray}}
\newcommand{\ea}{\end{eqnarray}}
\newcommand{\bas}{\begin{eqnarray*}}
\newcommand{\eas}{\end{eqnarray*}}
\newenvironment{pf}{\noindent{\sc Proof}.\enspace}{\rule{2mm}{2mm}\medskip}
\newenvironment{pfn}{\noindent{\sc Proof}}{\rule{2mm}{2mm}\medskip}

\newcommand{\R}{\mathbb{R}}

\author{ Cheikh Birahim NDIAYE$^{a}$, Abdul-Malik SAIID$^{b}$}

\date{}

\title{\bf An optimal boundary control approach to the Cherrier-Escobar problem}
\begin{document}

\newtheorem{lem}{Lemma}[section]
\newtheorem{pro}[lem]{Proposition}
\newtheorem{thm}[lem]{Theorem}
\newtheorem{rem}[lem]{Remark}
\newtheorem{cor}[lem]{Corollary}
\newtheorem{df}[lem]{Definition}

\maketitle

\begin{center}
{\small

\noindent  $^{a,b}$ Department of Mathematics Howard University \\  Annex $3$, Graduate School of Arts and Sciences, $217$\\ DC 20059 Washington, USA

}

\end{center}

\footnotetext[1]{E-mail addresses: cheikh.ndiaye@howard.edu, abdulmalik.saiid@bison.howard.edu\\
\thanks{\\ The first author was partially supported by NSF grant DMS--2000164.}}

\

\

\begin{center}
{\bf Abstract}
\end{center}
We study an optimal boundary control problem associated to the boundary obstacle problem for the couple conformal Laplacian and conformal Robin operator on \;$n$-dimensional compact Riemannian manifolds with boundary and with $n\geq 3$. When the Cherrier\cite{che1}-Escobar\cite{esc1} invariant of the compact Riemannian manifold with boundary is positive, we show that the optimal controls are equal to their associated optimal states. Moreover, we show that the optimal controls are minimizers of the Cherrier\cite{che1}-Escobar\cite{esc1} functional, and hence induce conformal metrics with zero scalar curvature and constant mean curvature. Furthermore, we show the existence of an optimal control under an Aubin\cite{aubin} type assumption. For the standard unit ball, we derive a sharp Sobolev trace type inequality and prove that the standard bubbles-namely conformal factor of metrics conformal to the standard one with zero scalar curvature and constant mean curvature- are the only optimal controls and hence equal to their associated optimal states.

 \begin{center}


\bigskip

\centerline{\bf AMS subject classification: 53C21, 35C60, 58J60, 55N10.}

\end{center}


\section{Introduction and statement of the results}\label{intro}
After the resolution of the Yamabe problem by Yamabe\cite{yamabe}, Trudinger\cite{tru}, Aubin\cite{aubin1}, and Schoen\cite{schoen}, Cherrier\cite{che1} and Escobar\cite{esc1} initiated an analogue of the Riemann mapping problem for manifolds with boundary that is referred to as the Cherrier-Escobar Problem. The Cherrier-Escobar problem is the geometric question which asks whether every compact Riemannian manifold with boundary of dimension greater than $2$ carries a conformal metric with zero scalar curvature and constant mean curvature. Analytically, it is equivalent to finding a smooth and positive solution to the geometric Boundary Value Problem (BVP)
\begin{equation}\label{eq:lequation}
\begin{cases}
L_gu=0 & \text{on }\; M, \\
B_gu=cu^{\frac{n}{n-2}} & \text{on }\; \partial M,
\end{cases}
\end{equation}
where $c$ is a real number, $(\overline{M}, g)$ is the background compact Riemannian manifold of dimension $n\geq3$\; with boundary \(\partial{M} \text{ and interior } M\), 
\[L_g:=-4\left(\frac{n-1}{n-2}\right)\D_g+R_g\] 
is the conformal Laplacian  of \( (\overline{M},g) \) and 
\[ B_g:= 4\left(\frac{n-1}{n-2}\right)\frac{\partial}{\partial n_g} + 2(n-1)H_g \] 
is the conformal Robin operator of \( (\overline{M},g) \), with \( R_g \) denoting the scalar curvature of \( (\overline{M},g) \), \( \Delta_g \) denoting the Laplace-Beltrami operator with respect to \( g \), \( H_g \) is the mean curvature of \( \partial M \) in \( (\overline{M},g) \), namely 
\[H_g := \frac{1}{n-1} tr_{\hat{g}} A_g \] 
where \( A_g \) is the second fundamental form of \( \partial M \) in \( (\overline{M}, g) \) with respect to the inner normal direction, \( \hat{g} := g|_{\partial M}, \, \frac{\partial}{\partial n_g} \) is the outer Neumann operator on \( \partial M \) with respect to \( g \). We set 
$$2^\#:=\frac{2(n-1)}{n-2}.$$
\vspace{5pt}

\noindent
The solution of \eqref{eq:lequation} given by the works of Cherrier\cite{che1}, Escobar\cite{esc1}, Marquez\cite{mar1}, Chen\cite{chen1}, and Almaraz\cite{alm1} 
is obtained by finding a smooth minimizer of the Cherrier\cite{che1}-Escobar\cite{esc1} functional $\mathcal{E}^g$ defined by
$$
\mathcal{E}^g(u):=\frac{\left<u, u\right>_g}{||u||^2_{L^{2^\#}(\partial M, \hat{g})}}, \;\;\;\;u\in H^1_+(\overline{M}, g):=\{u\in H^1(\overline{M}, g),\;u\geq0 \text{ in } M, \;u>0 \text{ on } \partial{M}\},
$$
where 
$$
\left<u,v\right>_{g}=4\frac{n-1}{n-2}\int_M\n_gu\n_gvdV_g+\int_{M}R_guvdV_g+2(n-1)\int_{\partial M}H_g uv dS_g, \;\;u, v\in H^1(\overline{M}, g),
$$
$dV_g$ is the Riemannian measure associated to \;$(\overline{M}, g)$, $dS_g$ is the Riemannian measure associated to $ (\partial M, \hat{g}), \; L^{2^\#}(\partial M, \hat{g})$ is the standard Lebesgue space of functions which are \;$2^\#$-integrable over \;$\partial{M}$\; with respect to \;$\hat{g}$, $||\cdot||_{L^{2^\#}(\partial M, \hat{g})}$\; is the standard  \;$L^{2^{\#}}$-norm on $L^{2^{\#}}(\partial M, \hat{g})$\; and \;$H^1(\overline{M}, g)$\; is the space of functions on \;$\overline{M}$\; which are \;$L^2$-integrable\; over $\overline{M}$ together with their first derivatives with respect to $g$, see \cite{aubin} for precise definitions.
We remark that for \;$u$\; smooth,
 $$
 \left<u,u\right>_{g}=\left<L_gu, u\right>_{L^2(\overline{M}, g)}+\left<B_gu, u\right>_{L^{2}(\partial M, \hat{g})},
 $$
where \;$\left<\cdot, \cdot\right>_{L^2(\overline{M}, g)}$\;denotes the $L^2$-scalar product on \;$\overline{M}$\; with respect to \;$g$ and \;$\left<\cdot, \cdot\right>_{L^2(\partial M, \hat{g})}$\;denotes the \;$L^2$-scalar product on \;$(\partial M, \hat{g})$\; with respect to \;$\hat{g}$.
\vspace{8pt}

\noindent
The Cherrier-Escobar  problem \eqref{eq:lequation} and the prescribed mean curvature problem ($K$ a smooth function) 
\begin{equation}\label{eq:sequation}
\begin{cases}
L_gu=0 & \text{on}\;\;M, \\
B_gu=Ku^{\frac{n}{n-2}} & \text{on}\;\;\partial M,
\end{cases}
\end{equation}
have been intensively studied using methods of Calculus of Variations, Critical Points Theory, Morse Theory, Dynamical Systems, Blow-up Analysis, Perturbation Methods, and Algebraic Topology, see \cite{alm1}, \cite{alm2}, \cite{aym}, \cite{chen1}, \cite{che1}, \cite{esc1}, \cite{ahm}, \cite{mar1}, and the references therein.
\vspace{10pt}

\noindent
In this paper, we investigate equation \eqref{eq:lequation} in the context of Optimal Control Theory. We recall that in the works of Cherrier\cite{che1}, Escobar\cite{esc1}, Marquez\cite{mar1}, Chen\cite{chen1}, and Almaraz\cite{alm1}, the most difficult case is the positive one, namely \;$\mathcal{\mu}(M,\partial M, [g])>0$, where
$$\mathcal{\mu}(M,\partial M, [g]):=\inf_{u\in H^1_+(\overline{M}, g) }\mathcal{E}^g(u)$$ is the Cherrier-Escobar invariant of $(M, \partial M, [g])$ for the zero scalar curvature and constant mean curvature problem and $$[g]=\{g_u=u^{\frac{4}{n-2}}g,\;\; u\in C^{\infty}_+(\overline{M}, g)\}$$  denotes the conformal class of $g$\;with
$$C^{\infty}_+(\overline{M}, g)=\{u\in C^{\infty}(\overline{M}, g):\;\;u>0\}$$ and \;$C^{\infty}(\overline{M}, g)$\; is the space of \;$C^{\infty}
$-functions on $\overline{M}$ with respect to $g$, see \cite{aubin}, \cite{lp} for a precise definition.
\vspace{7pt}

\noindent
In this paper, we focus on the case \;$\mathcal{\mu}(M, \partial M, [g])>0$. Precisely, under the assumption \;$\mathcal{\mu}(M, \partial M, [g])>0$, we study the following optimal control problem for the boundary obstacle problem associated to the couple $(L_g, B_g):$
\begin{equation}\label{stateop}
\text{Finding} \;\;\;u_{min}\in H^1_+(\overline{M}, g)\;\;\text{such that}\;\;\; \mathcal{I}^g(u_{min})=\min _{u\in H^1_+(\overline{M}, g)}\mathcal{I}^g(u),
\end{equation}
where
 $$
 \mathcal{I}^g(u)= \frac{ \left<u,u\right>_{g}}{||T_g(u)||_{L^{2^{\#}}(\partial M, \hat{g})}^2} ,\;\;\;\;u\in H^{1}_+(\overline{M}, g)
$$
with 
 \begin{equation}\label{defoo}
 T_g(u)=\arg{\min_{v\in H^1_+(\overline{M}, g),\\\;\;tr(v)\geq tr(u)}\left<v, v \right>_g}
\end{equation}
where the symbol $$\arg{\min_{v\in H^1_+(\overline{M}, g),\\\;\;trv\geq tru}\left<v, v \right>_g}$$ denotes the unique solution to the minimization problem \;$$\min_{v\in H^1_+(\overline{M}, g),\\\;\;trv\geq tru}\left<v, v \right>_g,$$ with $tr : H^1(\overline{M}, g) \to L^{2}(\partial{M}, \hat{g})$ denoting the trace operator, see Lemma \ref{obstacleq}.
\vspace{10pt}

\noindent
We obtain the following result for the optimal control problem \eqref{stateop}.

\begin{thm}\label{positivevarmass}
Assuming that  \;$\mathcal{\mu}(M, \partial M, \;[g])>0$, then\\
1) For $u\in H^1_+(\overline{M}, g)$, 
$$
\mathcal{I}^g(u)=\min _{v\in H^1_+(\overline{M}, g)}\mathcal{I}^g(v)\implies T_g(u)=u, \;\;u\in C^{\infty}_+(\overline{M}, g),\;\;\text{and}\;\;\begin{cases}
R_{g_u}=0\\
H_{g_u}=\text{const}>0,
\end{cases}
$$
with \;$g_u=u^{\frac{4}{n-2}}g$.\\\\
2) If there exists $ u_{min} \in  C^{\infty}_+(\overline{M}, g)$ such that $\mathcal{I}^g(u_{min})=\mathcal{\mu}_{oc}(M, \partial M, \;[g])$, then there exists also \;$u^{min}\in C^{\infty}_+(\overline{M}, g)$\; such that
$$
\mathcal{E}^g(u^{min})=\min _{v\in H^1_+(\overline{M}, g)}\mathcal{E}^g(v)\;\;\;{and}\;\;\ \begin{cases}
R_{g_{u^{min}}}=0\\
H_{g_{u^{min}}}=\mathcal{\mu}(M, \partial M, \;[g])\\
T_{g(u^{min})}=u^{min}
\end{cases}
$$
with \;$g_{u^{min}}=(u^{min})^{\frac{4}{n-2}}g$.\\\\
3) If
\begin{equation}
\mu_{oc}(M, \partial{M}, [g]) < \mu_{oc}(\mathbb{B}^n, \mathbb{S}^{n-1}, g_{\overline{\mathbb{B}}^n}) \tag{$a$}
\end{equation}
there exists $ u_{min} \in H_{+}^1(\overline{M}, g)$ such that $$\mathcal{I}^g(u_{min}) = \mu_{oc}(M, \partial{M}, [g])$$
\end{thm}
\vspace{7pt}

\noindent
 When \;$(\overline{M}, g)=(\overline{\mathbb{B}}^n, g_{\overline{\mathbb{B}}^n})$ is the standard closed \;$n$-dimensional unit ball of \;$\R^{n}$, we have the following result.
 
\begin{thm}\label{sphere}
Assuming that  \;$(\overline{M}, g)=(\overline{\mathbb{B}}^n, g_{\overline{\mathbb{B}}^n})$, then for $u\in H^1_+(\overline{M}, g)$,
$$
\mathcal{I}^g(u)=\min _{v\in H^1_+(\overline{M}, g)}I^g(v)\;\;\;\;\text{is equivalent to }\;\;\; u\in C^{\infty}_+(\overline{M}, g)\;\;\text{and}\;\; \begin{cases}
R_{g_u}=0\\
H_{g_u}=const>0,
\end{cases},
$$
with \;$g_u=u^{\frac{4}{n-2}}g$, \;$R_{g_{u}}$\;is the scalar curvature of \;$(\overline{M}, g_u)$, and $H_{g_u}$ is the mean curvature of \;$(\partial M, \hat{(g_u)})$.
\end{thm}
\begin{rem}\label{standardbubble}
We recall that the standard bubbles of the Cherrier-Escobar problem are the functions $u\in C^{\infty}_+(\overline{\mathbb{B}}^n, g_{\overline{\mathbb{B}}^n})$ with $R_{(g_{\overline{\mathbb{B}}^n})_u}=0$\;and $H_{(g_{\overline{\mathbb{B}}^n})_u}=\text{const}>0$\;where\;$(g_{\overline{\mathbb{B}}^n})_u=u^{\frac{4}{n-2}}g_{\overline{\mathbb{B}}^n}$. Thus, by point 1) of Theorem \ref{positivevarmass} and Theorem \ref{sphere}, we have $u$ is a standard bubble \;$\implies$\;$u=T_{g_{\overline{\mathbb{B}}^n}}(u)$.
\end{rem}
\begin{rem}\label{standardbubble}
It is well-known that the Cherrier–Escobar problem is intimately related to the $1/2$-Yamabe problem on conformal infinity of Poincaré–Einstein manifolds. We would like to emphasize that our methods also work for the $\gamma$-Yamabe problem (see \cite{ja}, \cite{aligo}, \cite{smw}, \cite{mayndia1}, \cite{nss}) on conformal infinity of asymptotically hyperbolic manifolds for arbitrary fractional parameter $\gamma$.
\end{rem}
\vspace{7pt}

\noindent
We will follow the same strategy as in the work of the second author \cite{nd8} to prove Theorem \ref{positivevarmass}-Theorem \ref{sphere}. We first use the variational characterization of the solution \;$T_g(u)$ of the boundary obstacle problem for the couple conformal Laplacian and conformal Robin operator (see Lemma \ref{obstacleq})\;  to show that \;$T_g$\; is  idempotent, see Proposition \ref{nilpotent}. Next, we use the idempotent property of \;$T_g$ to establish some monotonicity formulas for $\mathcal{E}^g$ and $\mathcal{I}^g$ involving $T_g$, see Lemma \ref{decreasingformula}, Lemma \ref{minimaltunnel}, and Lemma \ref{decreasingformulaop}. Using the latter monotonicity formulas, we show that any minimizer of \;$\mathcal{E}^g$\; and any solution of the optimal control problem \eqref{stateop} is a fixed point of \;$T_g$, see Corollary \ref{rigidityminimizer} and Corollary \ref{rigidityminimizeri}. After this, using the latter information, we derive that $$\inf_{H^1_+(\overline{M}, g)}\mathcal{I}^g=\inf_{H^1_+(\overline{M}, g)} \mathcal{E}^g $$ and
$$u\in H^1_+(\overline{M}, g) \text{ is a minimizer of } \;\mathcal{I}^g\; \Longleftrightarrow u\in H^1_+(\overline{M}, g) \text{ is a mimimizer of } \mathcal{E}^g,$$ see Proposition \ref{sameinf} and Proposition \ref{sameminimizer}. From this, Theorem \ref{positivevarmass} follows from an Aubin\cite{aubin} minimizing argument that is based on Ekeland's\cite{ivek1} variational principle and the analysis of Palais-Smale sequence of $\mathcal{E}^g$ by Almaraz \cite{alm3}. Furthermore, Theorem \ref{sphere} follows from the fact that the standard bubbles are the only minimizers of $\mathcal{E}^g$ on $H^1_+(\overline{M}, g)$ when $(\overline{M}, g)=(\overline{\mathbb{B}}^n, g_{\overline{\mathbb{B}}^n})$.
\vspace{7pt}

\noindent
The structure of the paper is as in \cite{nd8}. Indeed in Section \ref{eq:notpre}, we collect some preliminaries and fix some notations. In Section \ref{opp}, we discuss the boundary obstacle problem for the couple conformal Laplacian and conformal Robin operator. Moreover, we show that \;$T_g$\; is idempotent and homogeneous, and $\mathcal{I}^g$ is scale invariant. In Section \ref{conformalrule}, we derive the conformal transformation formulas of  \;$T_g$\;and \;$\mathcal{I}^g$. In Section \ref{monotonicity}, we show a monotonicity formula for $\mathcal{E}^g$ when passing from \;$u$\; to \;$T_g(u)$\; and discuss some applications. In Section \ref{monotonicityop}, we compare $\mathcal{I}^g$ and $\mathcal{E}^g$ and establish a monotonicity formula for $\mathcal{I}^g$ when passing from \;$u$\; to \;$T_g(u)$, and present some applications. In Section \ref{comparingval}, we show $$\inf_{H^1_+(\overline{M}, g)}\mathcal{I}^g=\inf_{H^1_+(\overline{M}, g)} \mathcal{E}^g$$ and \;$$u\in H^1_+(\overline{M}, g) \text{ is a minimizer of } \;\mathcal{I}^g\; \Longleftrightarrow \;u\; \text{ is a minimizer of } \mathcal{E}^g.$$ Moreover, we present the proof of Theorem \ref{positivevarmass}. Finally, in Section \ref{casesphere}, we derive a sharp Sobolev trace type inequality involving the obstacle operator  \;$T_g$\; and present the proof of Theorem \ref{sphere}.

\section{Notations and Preliminaries}\label{eq:notpre}
 In this brief section, we fix our notations and give some preliminaries. First of all, from now until the end of the paper, \;$(\overline{M}, g)$\; is the background compact Riemannian manifold with boundary of dimension $n\geq 3$. Even if Theorem \ref{positivevarmass} and Theorem \ref{sphere} are stated with the background metric $g$, we will be working, as in the work \cite{nd8}, with an arbitrary $h\in [g]$. 
\vspace{7pt}
 
\noindent
For  \;$h\in [g]$, we recall that the Cherrier-Escobar functional \;$\mathcal{E}^h$\; and its subcritical approximation \;$\mathcal{E}_p^h$ ($1\leq p<2^\#-1$) are defined by
 \begin{equation}\label{eq:defj1}
\mathcal{E}^h(u):= \frac{\left<u,u\right>_h}{||u||^2_{L^{2^\#}(\partial{M}, \hat{h})}} ,\;\;\;\;u\in H^1_+(\overline{M}, h),
\end{equation}
and
\begin{equation}\label{eq:defjt1}
\mathcal{E}^h_p(u):= \frac{\left<u,u\right>_h}{||u||^2_{L^{p+1}(\partial{M}, \hat{h})}} \;\;\;\;u\in H^1_+(\overline{M}, h),
\end{equation}
where
\;$H^1_+(\overline{M}, h)$, $\left<\cdot,\cdot\right>_h$, and  \;$L^{2^\#}(\partial{M}, \hat{h})$\; are as in Section \ref{intro} with \;$g$\; replaced by \;$h$ and \;$\hat{g}$\; replaced by \;$\hat{h}$. Moreover, $L^{p+1}(\partial{M}, \hat{h})$ is the standard Lebesgue space of functions which are $(p+1)$-integrable on \;$\partial{M}$\; with respect to \;$\hat{h}$ .
We set
\begin{equation}\label{valends}
\mathcal{E}_{2^\#-1}^h:=\mathcal{E}^h.
\end{equation}
and have a family of functional \;$(\mathcal{E}_p^h)_{1\leq p\leq 2^\#-1}$ defined on \;$H^1_+(\overline{M}, h)$. Clearly by definition of \;$\mathcal{E}^h_p$, 
\begin{equation}\label{phjph}
\mathcal{E}_p^h(\l u)=\mathcal{E}_p^h(u),\;\; \l>0,\; u\in H^1_+(\overline{M}, h).
\end{equation}
We recall that the conformal class $[g]$  of $g$ is
$$
[g]:=\{g_w:=w^{\frac{4}{n-2}}g, \;\;w\in C^{\infty}_+(\overline{M}, g)\},
$$
and in this paper we use the notation $$g_{*}:=*^{\frac{4}{n-2}}g.$$
\vspace{7pt}

\noindent
The following transformation rules are well-known and easy to verify (for \;$w\in C^{\infty}_+(\overline{M}, g)$)
\begin{equation}\label{inspace}
wC^{\infty}_+(\overline{M}, g_w)=C^{\infty}_+(\overline{M}, g), 
\end{equation}
\begin{equation}\label{inspace}
wH^{1}_+(\overline{M}, g_w)=H^1_+(\overline{M}, g).
\end{equation}
\begin{equation}\label{invv}
dV_{g_w}=w^{2^*}dV_g, \;\;dS_{g_w}=w^{2^{\#}}dS_g, \;\;2^*:=\frac{2n}{n-2}.
\end{equation}
\begin{equation}\label{invs}
\left<u, u\right>_{g_w}=\left<wu, wu\right>_{g}, \;\;\;u\in H^{1}_+(\overline{M}, g_w).
\end{equation}
\begin{equation}\label{invdf}
||u||_{L^{2^{\#}}(\partial{M}, g_w)}=||uw||_{L^{2^{\#}}(\partial{M}, g)}, \;\;\;u\in H^{1}_+(\overline{M}, g_w).
\end{equation}
and 
\begin{equation}\label{invyamabe}
\mathcal{E}^{g_w}(u)=\mathcal{E}^g(wu), \;\;\;u\in H^{1}_+(\overline{M}, g_w).
\end{equation}
We recall also the fact that \;$\inf_{H^1_+(\overline{M}, g) }\mathcal{E}^g$\; in the definition of \;$\mathcal{\mu}(M,\partial{M}, [g])$\; depends only on \;$[g]$\; can be seen from \eqref{inspace} and  \eqref{invyamabe}.
\vspace{7pt}

\noindent
For  \;$\mathcal{\mu}(M,\partial{M}, [g])>0$ and $h\in[g]$, we define the Cherrier-Escobar optimal boundary control functional $\mathcal{I}^h$ and its subcritical approximations\;$\mathcal{I}^h_p$ ($1\leq p<2^\#-1$) by 
\begin{equation}\label{eq:defi1}
\mathcal{I}^h(u):= \frac{\left<u,u\right>_h}{||T_h(u)||^2_{L^{2^\#}(\partial{M}, \hat{h})}} ,\;\;\;\;u\in H^1_+(\overline{M}, h),
\end{equation}
and
\begin{equation}\label{eq:defit1}
\mathcal{I}^h_p(u):= \frac{\left<u,u\right>_h}{||T_h(u)||^2_{L^{p+1}(\partial{M}, \hat{h})}}, \;\;\;\;u\in H^1_+(\overline{M}, h),
\end{equation}
where  \;$T_h$\; is as in \eqref{defoo} with \;$g$\; replaced by \;$h$.
We set
\begin{equation}\label{defend}
\mathcal{I}_{2^\#-1}^h:=\mathcal{I}^h
\end{equation}
and have a family of functional \;$(\mathcal{I}_p^h)_{1\leq p\leq 2^\#-1}$\; defined on \;$H^1_+(\overline{M}, h)$. 
\vspace{7pt}

\noindent
When $\mathcal{\mu}(M, \partial{M}, [g])>0$, the following Sobolev trace type inequality holds by definition of $\mu(M, \partial{M}, [g])$.
\begin{lem}\label{sobolev1}
Assuming \;$\mathcal{\mu}(M, \partial{M}, [g])>0$ and $h\in[g]$, then for $u\in H^1_+(\overline{M}, h)$
$$
||u||_{L^{2^{\#}}(\partial{M}, \hat{h})}\leq \frac{1}{\sqrt{\mathcal{\mu}(M, \partial{M}, [g])}}||u||_h
$$
\end{lem}
\vspace{7pt}

\noindent
When  \;$(\overline{M}, g)=(\overline{\mathbb{B}}^n,  g_{\overline{\mathbb{B}}^n})$, we have the following well-known stronger version of the latter Sobolev trace inequality, see \cite{esc1}.
\begin{lem}\label{sobolev2}
Assuming \;$(\overline{M}, g)=(\mathbb{B}^n, g_{\mathbb{B}^n})$ and $h=g_w\in[g]$, then for  $u\in H^1_+(\overline{M}, h)$,
\begin{equation}\label{inequality0}
||u||_{L^{2^{\#}}(\partial{M}, \hat{h})}\leq \frac{1}{\sqrt{\mathcal{\mu}(\mathbb{B}^n, \mathbb{S}^{n-1}, [g_{\overline{\mathbb{B}}^n}])}}||u||_h,
\end{equation}
with equality holding if and only if $$u\in C^{\infty}_+(\overline{M}, h)\;\;\text{and}\;\; \begin{cases}
R_{g_wu}=0\\
H_{g_wu}=const>0,
\end{cases},$$
\end{lem}
\begin{rem}
We recall that the explicit value of  \;$\mathcal{\mu}(\mathbb{B}^n, \mathbb{S}^{n-1}, [g_{\overline{\mathbb{B}}^n}])$\; is well-known, see \cite{esc1}.
\end{rem}

\section{Boundary obstacle problem for the couple $(L_h, B_h)$. }\label{opp}
In this section, we study the boundary obstacle problem for the couple conformal Laplacian and conformal Robin operator $(L_h, B_h)$ with \;$h\in [g]$\ under the assumption $\mathcal{\mu}(M, \partial{M}, [g])>0$. Indeed in the same spirit as in \cite{nd8}, given \;$u\in H^1_+(\overline{M}, h)$, we look for a solution to the minimization problem  
\begin{equation}\label{obspan}
\min_{v\in H^1_+(\overline{M}, h),\\\;\;tr(v)\geq tr(u)}\left<v, v \right>_h.
\end{equation} 
We anticipate that the assumption $\mathcal{\mu}(M, \partial{M}, [g])>0$ guarantees that the minimization problem \eqref{obspan} is well posed and $$||v||_h:=\sqrt{\left<v, v\right>_{h}}$$ is a norm on $H^1(\overline{M}, g)$, since 
$(L_h, B_h) \geq 0$ and $ker(L_h, B_h)={0}$ in the sense that the quadratic form $\left< \cdot, \cdot\right>_h$  associated to $(L_h, B_h)$ verifies $$\left<u, u\right>_h \geq 0, \; \forall u \in H^1(\overline{M}, g)$$ and $$ker(L_h, B_h):=\{ u \in H^1(\overline{M}, g) \; | \; \left<u, v\right>_h = 0, \; \forall v \in H^1(\overline{M}, g) \} = \{0\}$$ 
\vspace{7pt}

\noindent
We start with the following lemma providing the existence and uniqueness of solution for the boundary obstacle problem for the couple conformal Laplacian and conformal Robin operator \;$(L_h, B_h)$\; \eqref{obspan}. 
\begin{lem}\label{obstacleq}
Assuming that \;$\mathcal{\mu}(M, \partial{M}, [g])>0$ and $h\in[g]$, then  for \;$u\in H^1_+(\overline{M}, h)$,
there exists a unique $T_h(u)\in H^1_+(\overline{M}, h)$ such that
\begin{equation}\label{tg}
||T_h(u)||^2_h=\min_{v\in H^1_+(\overline{M}, h),\\\;\;tr(v)\geq tr(u)}||v||^2_h
\end{equation}
\end{lem}
\begin{pf}
Since \;$\mathcal{\mu}(M, \partial{M}, [g])>0$, then $(L_h, B_h) \geq 0$ and $ker(L_h, B_h)={0}$. Thus \;$<\cdot, \cdot>_h$\; defines a scalar product on \;$H^1_+(\overline{M}, h)$\; inducing a norm $||\cdot||_h$ equivalent to the standard \;$H^1(\overline{ M}, h)$-norm on \;$H^1_+(\overline{M}, h)$. Hence, as in the classical obstacle problem for the Laplace-Beltrami operator on closed Riemannian manifold, the lemma follows from the Direct Methods in the Calculus of Variations, and the fact that for every $u\in H^1_+(\overline{ M}, h)$, the set $\{v\in H^1_+(\overline{ M}, h): tr(v)\geq tr(u) \text{ on } \partial{M}\}$ is closed under the $H^1(\overline{ M}, h)$-weak convergence.
\end{pf}
\vspace{7pt}

\noindent
We study now some properties of the boundary obstacle solution map (or state map) \;$T_h\; : H^1_+(\overline{M}, h)\longrightarrow  H^1_+(\overline{M}, h)$. We start with the following algebraic one.
\begin{pro}\label{nilpotent}
Assuming that  \;$\mathcal{\mu}(M, \partial{M}, [g])>0$\; and  \;$h\in[g]$, then the state map \;$T_h\; : H^1_+(\overline{M}, h)\longrightarrow  H^1_+(\overline{M}, h)$\; is idempotent, i.e $$T^2_h=T_h.$$
\end{pro}

\begin{pf}
Let $v\in H^1_+(\overline{M}, h)$ such that $tr(v)\geq tr(T_h(u))$. Then $T_h(u)\geq u \text{ on } \partial{M}$ implies
$v\geq u \text{ on } \partial{M}$. Thus by minimality, \;$$||v||_h \geq ||T_h(u)||_h.$$ Since  \;$T_h(u)\in H^1_+(\overline{ M}, h)$\; and \;$T_h(u)\geq T_h(u)$\; on $\partial{M}$, then by uniqueness, we obtain that 
$$
T_h(T_h(u))=T_h(u),
$$
which ends the proof.
\end{pf}
\vspace{7pt}

\noindent
We have the following lemma showing that $T_h$ is positively homogeneous.
\begin{lem}\label{pht}
Assuming  \;$\mathcal{\mu}(M, \partial{M}, [g])>0$\; and  \;$h\in [g]$, then for  \;$\l>0$, 
$$
T_h(\lambda u)=\lambda T_h(u), \;\;\;\forall u\in H^1_+(\overline{M}, h).
$$
\end{lem}
\begin{pf}
Let $v \in H^1_+(\overline{M}, h) \text{ such that } \l u \leq v \text{ on } \partial{M}$. Then $\l>0$ implies $\l^{-1}v\geq u$. Thus, since \;$\l^{-1}v\in H^1_+(\overline{M}, h)$, then again by minimality 
$$
||\l^{-1}v||_h\geq ||T_h(u)||_h.
$$
Hence by positive homogeneity of $||\cdot||_h$, we get
$$
||v||_h\geq ||\l T_h(u)||_h
$$
And, since $\l T_h(u)\in H^1_+(\overline{M}, h)$ and $\l T_h(u) \geq \l u \; \text{on} \; \partial{M}$ , then by uniqueness we obtain
$$
T_h(\l u)=\l T_h(u),
$$
as desired.
\end{pf}
\vspace{7pt}

\noindent
Lemma \ref{pht} implies the following analogue of formula \eqref{phjph} for $\mathcal{I}^h_p$.
\begin{cor}\label{scaleinv}
Assuming  \;$\mathcal{\mu}(M, \partial{M}, [g])>0$,  \;$h\in [g]$, and \;$1\leq p\leq 2^\#-1$, then for \;$\l>0$,
$$
\mathcal{I}^h_p(\lambda u)=\mathcal{I}_p^h(u), \;\;\;\forall u\in H^1_+(\overline{ M}, h).
$$
\end{cor}
\begin{pf}
It follows directly from the definition of \;$\mathcal{I}^h_p$ (see \eqref{eq:defi1}-\eqref{defend}), Lemma \ref{pht}), and the positive homogeneity of norms. Indeed, \eqref{eq:defi1}-\eqref{defend} imply
\begin{equation}\label{invs1}
\mathcal{I}^h_p(\lambda u)=\frac{||\l u||^2_h}{||T_h(\l u)||^2_{L^{p+1}(\partial{M}, \hat{h})}}. 
\end{equation}
Thus, Lemma \ref{pht} and \eqref{invs1} give
\begin{equation}\label{invs2}
\mathcal{I}^h_p(\lambda u)=\frac{||\l u||^2_h}{||\l T_h(u)||^2_{L^{p+1}(\partial{M}, \hat{h})}}.
\end{equation}
So, the positive homogeneity of norms and \eqref{invs2} imply
\begin{equation}\label{invs3}
\mathcal{I}^h_p(\lambda u)=\frac{||u||^2_h}{||T_h(u)||^2_{L^{p+1}(\partial{M}, \hat{h})}} .
\end{equation}
Hence, using again \eqref{eq:defi1}-\eqref{defend} combined with \eqref{invs3}, we get
$$
\mathcal{I}^h_p(\lambda u)=\mathcal{I}^h_p(u),
$$
as desired.
\end{pf}
\section{Conformal transformation rules of  \;$T_{h}$\; and  \;$\mathcal{I}^{h}$.}\label{conformalrule}
Here, we study the conformal transformation rules of \;$T_h$\; when \;$h$\; varies in \;$[g]$\; and use it to establish an analogue of formula \eqref{invyamabe} for the Cherrier-Escobar optimal control function $\mathcal{I}^h$. 
\vspace{7pt}

\noindent
We adopt the notation \;$h=g_w=w^{\frac{4}{n-2}}g$\; with \;$w\in C^{\infty}_+(\overline{M}, g)$\; and have:
\begin{lem}\label{tth}
Assuming that  \;$\mathcal{\mu}(M, \partial{M}, [g])>0$\; and \;$w\in C^{\infty}_+(\overline{M}, g)$, then
$$
T_{g_w}(u)=w^{-1}T_g(wu), \;\;\;u\in H^1_+(\overline{M}, g_w).
$$
\end{lem}
\begin{pf}
The proof here follows the same strategy employed in \cite{nd8}. Indeed, let \;$v\in H^1_+(\overline{M}, g_w)$\; with \;$tr(v)\geq tr(u)$. Then using formula \eqref{invs}, we get
$$
||v||_{g_w}=||wv||_g.
$$
Hence, since \;$tr(v)\geq tr(u)$\; and \;$w>0$\; implies \;$tr(vw)\geq tr(uw)$, then  using $wu, \;wv\in wH^1_+(\overline{M}, g_w)$ and  \eqref{inspace}, by minimality we have
$$
||v||_{g_w}\geq ||T_g(wu)||_{g}.
$$
So, using again formula \eqref{invs}, we obtain
$$
||v||_{g_w}\geq ||w^{-1}T_g(wu)||_{g_w}
$$
Now, since  $w^{-1}T_g(wu)\geq w^{-1}wu=u \text{ on } \partial{M}$ and $w^{-1}T_g(wu)\in H^1_+(\overline{M}, g_w)$ (see \eqref{inspace}), then by uniqueness 
$$
T_{g_w}(u)=w^{-1}T_g(wu),
$$ 
thereby ending the proof.
\end{pf}
\vspace{7pt}

\noindent
As a consequence of Lemma \ref{tth}, we have:
\begin{cor}\label{invfixpt}
Assuming that  \;$\mathcal{\mu}(M, \partial{M}, [g])>0$\; and  \;$w\in C^{\infty}_+(\overline{M}, g)$, then
$$
Fix(T_{g_w})=w^{-1}Fix(T_g),
$$
where for $h\in [g]$, and 
\begin{equation}\label{deffix}
Fix(T_h):=\{u\in H_+^1(\overline{M}, h):\;\;T_h(u)=u\}.
\end{equation}
\end{cor}
\begin{pf}
Using lemma \ref{tth} and \eqref{deffix}, we have
 $$
u\in Fix(T_{g_w})\Longleftrightarrow\;\;T_{g_w}(u)=u\;\;\Longleftrightarrow\;\; w^{-1}T_g(wu)=u.
$$ 
Thus
$$u\in Fix(T_{g_w})\;\;\Longleftrightarrow\;\;T_g(wu)=wu\;\;\Longleftrightarrow \;\;wu\in Fix(T_g).$$
So$$
u\in Fix(T_{g_w})\;\;\;\Longleftrightarrow\;\;\;u\in w^{-1}Fix(T_g),
$$
as desired
\end{pf}
\vspace{7pt}

\noindent
Lemma \ref{tth} also implies the following corollary.
\begin{cor}\label{invd}
Assuming that  \;$\mathcal{\mu}(M, \partial{M}, [g])>0$\; and \;$w\in C^{\infty}_+(\overline{M}, g)$, then
$$
||T_{g_w}(u)||_{L^{2^\#}(\partial{M}, \hat{(g_w)})}=||T_{g}(wu)||_{L^{2^\#}(\partial{M}, \hat{g})}, \;\;\;u\in H^1_+(\overline{M}, g_w).
$$
\end{cor}
\begin{pf}
Using the definition of \;$L^{2^\#}(\partial{M}, \hat{(g_w)})$-norm, we have
$$
||T_{g_w}(u)||_{L^{2^\#}(\partial{M}, \hat{(g_w)})}^{2^\#}=\int_{\partial{M}}T_{g_w}(u)^{2^\#}dS_{g_w}.
$$
Thus, \eqref{invv} and Lemma \ref{tth}, imply
$$
||T_{g_w}(u)||_{L^{2^\#}(\partial{M}, \hat{(g_w)})}^{2^\#}=\int_{\partial{M}}T_{g}(wu)^{2^\#}dS_{g_w}
$$
So, by the definition of the \;$L^{2^\#}(M, g)$-norm, we obtain
$$
||T_{g_w}(u)||_{L^{2^\#}(\partial{M}, \hat{(g_w)})}^{2^\#}=||T_{g}(wu)||_{L^{2^\#}(\partial{M}, \hat{g})}^{2^\#}.
$$
Hence, we get 
$$
||T_{g_w}(u)||_{L^{2^\#}(\partial{M}, \hat{(g_w)})}=||T_{g}(wu)||_{L^{2^\#}(\partial{M}, \hat{g})},
$$
as desired.
\end{pf}
\vspace{7pt}

\noindent
As a consequence of Corollary \ref{invd}, we have the following analogue of formula \eqref{invyamabe} for the Cherrier-Escobar optimal boundary control functional.
\begin{cor}\label{invi}
Assuming that  \;$\mathcal{\mu}(M, \partial{M},  [g])>0$\; and \;$w\in C^{\infty}_+(\overline{M}, g)$, then for  \;\;$u\in H^1_+(\overline{M}, g_w)$, 
$$
\mathcal{I}^{g_w}(u)=\mathcal{I}^g(wu).
$$
\end{cor}
\begin{pf}
Using definition of \;$\mathcal{I}^{g_w}$ (see \eqref{eq:defi1}), we have
$$
\mathcal{I}^{g_w}(u)=\frac{||u||^2_{g_w}}{||T_{g_w}(u)||^2_{L^{2^\#}(\partial{M}, \hat{(g_w)})}}
$$
Thus, \eqref{invs} and Corollary \ref{invd}, imply
$$
\mathcal{I}^{g_w}(u)=\frac{||wu||^2_{g}}{||T_g(wu)||^2_{L^{2^\#}(\partial{M}, \hat{(g_w)})}}
$$
So by the definition of \;$\mathcal{I}^g$ (see \eqref{eq:defi1}), we obtain
$$
\mathcal{I}^{g_w}(u)=\mathcal{I}^g(wu),
$$
as desired.	
\end{pf}
\vspace{5pt}

\noindent
From Corollary \ref{invi} we have.
\begin{cor}\label{cor14}
Assuming that  \;$\mathcal{\mu}(M, \partial{M}, [g])>0$\; and \;$w\in C^{\infty}_+(\overline{M}, g)$, then
$$
\inf_ {u\in H^1_+(\overline{M}, g_w)} \mathcal{I}^{g_w}(u)=\inf _ {u\in H^1_+(\overline{M}, g)} \mathcal{I}^g(u).
$$
\end{cor}
\begin{pf}
Using \eqref{inspace} we have
 $$
 \inf _ {u\in H^1_+(\overline{M}, g)} \mathcal{I}^g(u)=\inf _ {u/w\in H^1_+(\overline{M}, g_w)} \mathcal{I}^g(u)
 $$
 Thus, setting $u=w\bar u$, we obtain
 $$
 \inf _ {u\in H^1_+(\overline{M}, g)} \mathcal{I}^g(u)=\inf _ {\bar u\in H^1_+(\overline{M}, g_w)} \mathcal{I}^g(w\bar u)
 $$
 Hence, Corollary \ref{invi} implies
 $$
 \inf _ {u\in H^1_+(\overline{M}, g)} \mathcal{I}^g(u)=\inf _ {\bar u\in H^1_+(\overline{M}, g_w)} \mathcal{I}^{g_w}(\bar u),
 $$
as desired.
 \end{pf}
\vspace{7pt}

\noindent
Similar to the Cherrier-Escobar invariant \;$\mathcal{\mu}(M,\partial{M}, [g])$, Corollary \ref{cor14} justifies the following definition.
\begin{df}\label{defyp}
Assuming that  \;$\mathcal{\mu}(M, \partial{M}, [g])>0$ and \;$h\in [g]$, then
$$
\mathcal{\mu}_{oc}(M, \partial{M}, [h]):=\inf _ {u\in H^1_+(\overline{M}, h)} \mathcal{I}^h(u).
$$
\end{df}
\vspace{7pt}

\noindent
\begin{rem}
Clearly \;for $h\in [g]$, $$\mathcal{\mu}_{oc}(M,\partial{M}, [h])=\mathcal{\mu}_{oc}(M,\partial{M}, [g]).$$
\end{rem}


\vspace{7pt}

\noindent
In the works of Cherrier\cite{che1} and Escobar\cite{esc1}, the following family of real numbers (not conformal invariant) was introduced
$$
\mathcal{\mu}^p(\overline{M}, h):=\inf _ {u\in H^1_+(\overline{M}, h)} \mathcal{E}^h_p(u), \;\;h\in[g], \;1\leq p<2^\#-1.
$$
We set
$$
\mathcal{\mu}^{2^\#-1}(M, \partial{M}, h)=\mathcal{\mu}(M, \partial{M}, [h])
$$
and have a family of real numbers $(\mathcal{\mu}^p(M, \partial{M}, h))_{1\leq p\leq 2^\#-1}$ which is conformally invariant for $p=2^\#-1$.
\vspace{7pt}

\noindent
Similarly, for $\mathcal{\mu}(M, \partial{M}, [g])>0$ we define
$$
\mathcal{\mu}^p_{oc}(M, \partial{M}, h):=\inf _ {u\in H^1_+(\overline{M}, h)} \mathcal{I}^h_p(u), \;\;h\in[g], \;1\leq p<2^\#-1
$$
and
$$
\mathcal{\mu}^{2^\#-1}_{oc}(M, \partial{M}, h)=\mathcal{\mu}_{oc}(M, \partial{M}, [h]).
$$
This defines a family of real numbers $(\mathcal{\mu}^p_{oc}(M, \partial{M}, h))_{1\leq p\leq 2^\#-1}$ which is conformally invariant for $p=2^\#-1$
\section{Monotonicity formula for \;$\mathcal{E}^h_p$, $h\in [g]$, and $1\leq p\leq 2^\#-1$.}\label{monotonicity}
In this section, we present a monotonicity formula for \;$\mathcal{E}^h_p$\; when passing from \;$u$\; to \;$T_h(u)$\; for \;$h\in[g]$\; with \;$\mathcal{\mu}(M, \partial{M}, [g])>0$ \; and \;$1\leq p\leq 2^\#-1$. Moreover, we present some applications on the relation between the ground state of \;$\mathcal{E}^h_p$\; and the the fixed points of \;$T_h$.  
\vspace{7pt}

\noindent
The monotonicity formula reads as follows.
\begin{lem}\label{decreasingformula}
Assuming that  \;$\mathcal{\mu}(M, \partial{M}, [g])>0$, \;$h\in[g]$\; and \;$1\leq p\leq 2^{\#}-1$, then for \;$u\in H^1_+(\overline{M}, h)$, $$
\mathcal{E}_p^h(u)-\mathcal{E}_p^h(T_h(u))\geq \frac{1}{||T_h(u)||^2_{L^{p+1}(\partial{M}, \hat{h})}}\left[||u||^2_h-||T_h(u))||^2_h\right]\geq 0.
$$
\end{lem}
\begin{pf}
From the definition of \;$\mathcal{E}_p^h$ (see \eqref{eq:defj1}-\eqref{valends}), we have
\begin{equation}
\mathcal{E}_p^h(u)-\mathcal{E}_p^h(T_h(u))= \frac{||u||^2_h}{||u||^2_{L^{p+1}(\partial{M}, \hat{h})}}-\frac{||T_h(u)||^2_h}{||T_h(u)||^2_{L^{p+1}(\partial{M}, \hat{h})}}.
\end{equation}
Hence the result follows from  \;$tr(T_h(u))\geq tr(u)> 0$, and and the definition of \;$T_h$\; (see Lemma \ref{obstacleq}).
\end{pf}
\vspace{7pt}

\noindent
Lemma \ref{decreasingformula} implies the following rigidity result.
\begin{cor}\label{rigidity}
Assuming that  \;$\mathcal{\mu}(M, \partial{M}, [g])>0$, \;$h\in[g]$, and \;$1\leq p\leq 2^{\#}-1$, then for\;$u\in H^1_+(\overline{M}, h)$,
\begin{equation}\label{ineq}
\mathcal{E}_p^h(T_h(u)) \leq  \mathcal{E}_p^h(u)
\end{equation}
and  
\begin{equation}\label{eql}
\mathcal{E}_p^h(u)=\mathcal{E}_p^h(T_h(u))\implies u \in Fix(T_h).
\end{equation}
\end{cor}
\begin{pf}
Using Lemma \ref{decreasingformula}, we have
\begin{equation}\label{eq1}
\mathcal{E}_p^h(u)-\mathcal{E}_p^h(T_h(u))\geq  \frac{1}{||T_h(u)||^2_{p+1}}\left[||u||^2_h-||T_h(u)||^2_h\right]\geq 0.
\end{equation}
Thus, \eqref{ineq} follows from \eqref{eq1}. If \;$\mathcal{E}_p^h(u)=\mathcal{E}_p^h(T_h(u))$, then \eqref{eq1} implies
$$
||u||^2_h=||T_h(u)||^2_h.
$$
Hence, since \;$tr(u)\geq tr(u) \text{ and } u \in H^1_+(\overline{M}, h)$, then the uniqueness part in Lemma \ref{obstacleq} implies
$$
u=T_h(u).
$$
Hence, using \eqref{deffix}, we have
$$
u\in Fix(T_h),
$$
thereby ending the proof of the corollary.
\end{pf}


\vspace{7pt}

\noindent
Corollary \ref{rigidity} implies that minimizers of  \;$\mathcal{E}_p^h$\; on \;$H^1_+(\overline{M}, h)$\; belongs to \;$Fix(T_h)$. Indeed, we have:
\begin{cor}\label{rigidityminimizer}
Assuming that  \;$\mathcal{\mu}(M, \partial{M}, [g])>0$,  $h\in[g]$, and  \;$1\leq p\leq 2^{\#}-1$, then for \;$u\in H^1_+(\overline{M}, h)$, $$\mathcal{E}_p^h(u)=\mathcal{\mu}^p(M, \partial{M}, h)\implies u\in Fix(T_h).$$
\end{cor}
\begin{pf}
$\mathcal{E}_p^h(u)=\mathcal{\mu}^{p}(M ,\partial{M}, h)$\; implies
\begin{equation}\label{eq3}
\mathcal{E}_p^h(u)\leq \mathcal{E}_p^h(T_h(u)).
\end{equation}
Thus combining \eqref{ineq} and \eqref{eq3}, we get
\begin{equation}\label{eq4}
\mathcal{E}_p^h(u)= \mathcal{E}_p^h(T_h(u)).
\end{equation}
Hence, combining \eqref{eql} and \eqref{eq4}, we obtain
$$
u\in Fix (T_h),
$$
as desired
\end{pf}

\vspace{7pt}

\noindent
\begin{rem}\label{rem2022}
Under the assumption of Corollary \ref{rigidity}, we have Proposition \ref{nilpotent} and Corollary \ref{rigidity} imply that we can assume without loss of generality that any minimizing sequence \;$(u_l)_{l\geq 1}$\; of \;$\mathcal{E}_p^h$\; on  \;$H^1_+(\overline{M}, h)$\;satisfies
$$
u_l\in Fix(T_h), \;\;\;\;\forall l\ge 1.
$$
Indeed, suppose \;$u_l$\; is a minimizing sequence for \;$\mathcal{E}_p^h$\; on \;$H^1_+(\overline{M}, h)$. Then \;$u_l\in H^1_+(\overline{M}, h)$\; and 
$$
\mathcal{E}_p^h(u_l)\longrightarrow \inf_{H^1_+(\overline{M}, h)}\mathcal{E}_p^h.
$$
Thus by definition of infimum and Corollary \ref{rigidity}, we  have
$$
\inf_{H^1_+(\overline{M}, h)}\mathcal{E}_p^h\leq \mathcal{E}_p^h(T_h(u_l))\leq \mathcal{E}_p^h(u_l).
$$
This implies $$\mathcal{E}_p^h(T_h(u_l))\longrightarrow \inf_{H^1_+(\overline{M}, h)}\mathcal{E}_p^h.$$ Hence setting $$\hat u_l=T_h(u_l),$$ and using Proposition \ref{nilpotent}, we get
$$\mathcal{E}_p^h(\hat u_l)\longrightarrow \inf_{H^1_+(\overline{M}, h)}\mathcal{E}_p^h\;\;\; \;\;\text{and}\;\; \;\;\;\;\hat u_l=T_h(\hat u_l)$$ as desired.
\end{rem}

\section{Monotonicity formula for $\mathcal{I}^h_p$, $h\in [g]$,  $1\leq p\leq 2^\#-1$.}\label{monotonicityop}
In this section, we derive a monotonicity formula for $\mathcal{I}_p^h$ similar to the one for \; $\mathcal{E}_p^h$ \; derived in the previous section for $h\in [g]$\; with \;$\mathcal{\mu}(M, \partial{M}, [g])>0$\; and \;$1\leq p\leq 2^{\#}-1$. Moreover, we present some applications for \;$\mathcal{I}^h_p$\; similar to the ones done for $\mathcal{E}^h_p$ in the previous section.
\vspace{7pt}

\noindent
We have the following monotonicity formula for the \;$p$-Cherrier-Escobar optimal boundary control functional \;$\mathcal{I}_p^h$.
\begin{lem}\label{decreasingformulaop}
Assuming that  \;$\mathcal{\mu}(M, \partial{M}, [g])>0$, $h\in [g]$, and  \;$1\leq p\leq 2^{\#}-1$, then for  \;$u\in H^1_+(\overline{M}, h)$,
$$
\mathcal{I}_p^h(u)-\mathcal{I}_p^h(T_h(u))=\frac{1}{||T_h(u)||^2_{L^{p+1}(\partial{M}, \hat h)}}\left[||u||^2_h-||T_h(u)||^2_h\right]\geq 0.
$$
\end{lem}
\begin{pf}
Using the definition of \;$\mathcal{I}_p^h$ (see \eqref{eq:defi1}-\eqref{defend}), we have
$$
\mathcal{I}_p^h(u)-\mathcal{I}_p^h(T_h(u))=\frac{||u||^2_h}{||T_h(u)||^2_{L^{p+1}(\partial{M}, \hat h)}}-\frac{||T_g(u)||^2_h}{||T_h^2(u)||^2_{L^{p+1}(\partial{M}, \hat h)}}.
$$
Using \;$T_h^2(u)=T_h(u)$\;(see Proposition \ref{nilpotent}) and the definition of \;$T_h$\; (see Lemma \ref{obstacleq}), we  get
$$\mathcal{I}_p^h(u)-\mathcal{I}_p^h(T_h(u))=\frac{1}{||T_h(u)||^2_{L^{p+1}(\partial{M}, \hat h)}}\left[||u||^2_h-||T_g(u)||^2_h\right]\geq 0,$$
which ends the proof.
\end{pf}
\vspace{7pt}

\noindent
Similar to the previous section,  Lemma \ref{decreasingformulaop} implies the following rigidity result.
\begin{cor}\label{rigidityop}
Assuming that \;$\mathcal{\mu}(M, \partial{M}, [g])>0$,  $h\in[g]$, and  \;$1\leq p\leq 2^{\#}-1$, then for  \;$u\in H^1_+(\overline{M}, h)$,
\begin{equation}\label{ineqop}
\mathcal{I}_p^h(T_h(u)) \leq  \mathcal{I}_p^h(u)
\end{equation}
and  
\begin{equation}\label{eqlop}
\mathcal{I}_p^h(u)=\mathcal{I}_p^h(T_h(u))\implies u \in Fix(T_h).
\end{equation}
\end{cor}
\begin{pf}
From Lemma \ref{decreasingformulaop}, we have
\begin{equation}\label{eq1op}
\mathcal{I}_p^h(u)-\mathcal{I}_p^h(T_h(u))=\frac{1}{||T_h(u)||^2_{p+1}}\left[||u||^2_h-||T_h(u)||^2_h\right]\geq 0.
\end{equation}
So, \eqref{ineqop} follows from \eqref{eq1op}. Observe that if \;$\mathcal{I}_p^h(u)=\mathcal{I}_p^h(T_h(u))$, then \eqref{eq1op} implies
$$
||u||^2_h=||T_h(u)||^2_h.
$$
Hence, since \;$u\in H^1_+(\overline{M}, h)$ and $tr(u)\geq tr(u)$ on $\partial{M}$, then  as in the previous section the uniqueness part in Lemma \ref{obstacleq} implies
$$
u\in Fix(T_h),
$$
which ends the proof of the corollary.
\end{pf}
\vspace{7pt}

\noindent
As in the previous section, corollary \ref{rigidity} imply that minimizers of  \;$\mathcal{I}_p^h$\; belongs to  \;$Fix(T_h)$.
\begin{cor}\label{rigidityminimizeri}
Assuming that  \;$\mathcal{\mu}(M, \partial{M}, [g])>0$,  $h\in[g]$\;  and  \;$1\leq p\leq 2^{\#}-1$, then for \;$u\in H^1_+(\overline{M}, h)$, $$\mathcal{I}_p^h(u)=\mathcal{\mu}_{oc}^p(M, \partial{M}, h)\implies u \in Fix(T_h).$$
\end{cor}
\begin{pf} \text{ The proof is similar to the one in \cite{nd8}. Indeed }
$\mathcal{I}_p^h(u)=\mathcal{\mu}^p_{oc}(M, \partial{M}, h)$\; implies
\begin{equation}\label{eq3op}
\mathcal{I}_p^h(u)\leq \mathcal{I}_p^h(T_h(u)).
\end{equation}
Thus combining \eqref{ineqop} and \eqref{eq3op}, we get
\begin{equation}\label{eq4op}
\mathcal{I}_p^h(u)= \mathcal{I}_p^h(T_h(u)).
\end{equation}
Hence, combining \eqref{eqlop} and \eqref{eq4op}, we obtain
$$
u\in Fix (T_h),
$$
which ends the proof.

\end{pf}
\vspace{7pt}

\noindent
\begin{rem}
As in \cite{nd8}, under the assumptions of Corollary \ref{decreasingformulaop} and using the same argument as in Remark \ref{rem2022}, we have that Proposition \ref{nilpotent} and Corollary \ref{decreasingformulaop} imply that for a minimizing sequence \;$(u_l)_{l\geq 1} $\; of \;$\mathcal{I}_p^h$\; on  \;$H^1_+(\overline{M}, h)$, we can  assume without loss of generality that
$$
u_l\in Fix(T_h), \;\;\;\;\forall l\ge 1.
$$
\end{rem}
\vspace{10pt}

\noindent
\section{Comparing \;$\mathcal{\mu}^p(M, \partial{M}, h)$\; and \;$\mathcal{\mu}_{oc}^p(M,\partial{M}, h)$, $h\in[g]$, and \;$1\leq p\leq 2^\#-1$}\label{comparingval}
Here, for \;$\mathcal{\mu}(M,\partial{M}, [g])>0$,  \;$h\in [g]$,  \;$1\leq p\leq 2^\#-1$, we show  \;$\mathcal{\mu}^p(M, \partial{M}, h)=\mathcal{\mu}_{oc}^p(M, \partial{M}, h)$\; and \;$\mathcal{E}^h_p(u)=\mathcal{\mu}^p(M, \partial{M}, h)\Longleftrightarrow \mathcal{I}^h_p(u)=\mathcal{\mu}_{oc}^p(M, \partial{M}, h)$. As a consequence, we deduce Theorem \ref{positivevarmass}. 
\vspace{7pt}

\noindent
We start with the following comparison result by showing that  \;$\mathcal{I}_p^h\leq \mathcal{E}_p^h$\; and \;$\mathcal{I}_p^h=\mathcal{E}_p^h$\; on \;$T_h(H^1_+(\overline{M}, h))$\; the range of \;$T_h$. 
\begin{lem}\label{minimaltunnel}
Assuming that \;$\mathcal{\mu}(M, \partial{M}, [g])>0$,  $h\in[g]$\;  and  \;$1\leq p\leq 2^{\#}-1$, then 
\begin{equation}\label{lesval}
\mathcal{I}_p^h\leq \mathcal{E}_p^h\;\;\;\;\;\;\text{on}\;\;\;\;\; H^1_+(\overline{M}, h)
\end{equation}
and
\begin{equation}\label{sameval}
\mathcal{E}_p^h\circ T_h=\mathcal{I}_p^h\circ T_h\;\;\;\;\text{on}\;\;\; H^1_+(\overline{M}, h).
\end{equation}
\end{lem}
\begin{pf}
Using the definition of \;$\mathcal{E}_p^h$\; and \;$\mathcal{I}_p^h$ (see \eqref{eq:defj1}-\eqref{valends}, and \eqref{eq:defi1})-\eqref{defend},we have, for $u\in H^1_+(\overline{M}, h)$ 
\begin{equation}\label{app1}
\mathcal{E}_p^h(u)-\mathcal{I}_p^h(u)=\frac{||u||^2_h}{||u||^2_{L^{p+1}(\partial{M}, \hat h)}}-\frac{||u||^2_h}{||T_h(u)||^2_{L^{p+1}(\partial{M}, \hat h)}}.
\end{equation}
Thus \eqref{lesval} follows from \;$tr(T_h(u))\geq tr(u)>0$ and \eqref{app1}. Furthermore, we have
$$
\mathcal{E}_p^h(T_h(u))-\mathcal{I}_p^h(T_h(u))= \frac{||T_h(u)||^2_h}{||T_h(u)||^2_{L^{p+1}(\partial{M}, \hat h)}}-\frac{||T_h(u)||^2_h}{||T_h^2(u)||^2_{L^{p+1}(\partial{M}, \hat h)}}.
$$
Hence, \;$T_h^2(u)=T_h(u)$\; (see Proposition \ref{nilpotent}) implies \;$$\mathcal{E}_p(T_h(u))=\mathcal{I}_p^h(T_h(u)),$$
which ends the proof.
\end{pf}
\begin{rem}
Clearly \;$tr(T_h(u))\geq tr(u)>0$ and \eqref{app1} imply
$$
\mathcal{I}_p^h(u)=\mathcal{E}_p^h(u)\Longleftrightarrow tr(u)=tr(T_h(u)) \quad \text{on } \partial{M}
$$
\end{rem}
\vspace{7pt}

\noindent
Corollary \ref{rigidity}, Corollary \ref{rigidityop}, and Lemma \ref{minimaltunnel} imply \;$\mathcal{\mu}^p(M, \partial{M}, h)=\mathcal{\mu}^p_{oc}(M, \partial{M}, h)$. Indeed, we have:
\begin{pro}\label{sameinf}
Assuming that \;$\mathcal{\mu}(M, \partial{M}, [g])>0$,  \;$h\in[g]$\;  and  \;$1\leq p\leq 2^{\#}-1$, then
$$
\mathcal{\mu}^p(M, \partial{M}, h)=\mathcal{\mu}^p_{oc}(M, \partial{M}, h).
$$
\end{pro}
\begin{pf}
From \eqref{ineq}, we obtain
\begin{equation}\label{relinf1}
\mathcal{\mu}^p(M, \partial{M}, h)=\inf_{H^1_+(\overline{M}, h)}\mathcal{E}_p^h\circ T_h.
\end{equation}
Similarly, \eqref{ineqop} implies
\begin{equation}\label{relinf2}
\mathcal{\mu}^p_{oc}(M, \partial{M}, h)=\inf_{H^1_+(\overline{M}, h)}\mathcal{I}_p^h\circ T_h.
\end{equation}
Therefore, the result follows from \eqref{sameval}, \eqref{relinf1}, and \eqref{relinf2}.
\end{pf}
\vspace{7pt}

\noindent
Combining Corollary \ref{rigidityminimizer}, Corollary \ref{rigidityminimizeri}, \eqref{sameval} and Proposition \ref{sameinf}, we have.
\begin{pro}\label{sameminimizer}
Assuming that \;$\mathcal{\mu}(M,\partial{M}, [g])>0$,  $h\in[g]$\;  and  \;$1\leq p\leq 2^{\#}-1$, then for \;$u\in H^1_+(\overline{M}, h)$, 
$$
\mathcal{E}^h_p(u)=\mathcal{\mu}^p(M, \partial{M}, h)\;\;\;\Longleftrightarrow\;\; \mathcal{I}^h_p(u)=\mathcal{\mu}^p_{oc}(M, \partial{M}, h).
$$
\end{pro}
\begin{pf}
If $\mathcal{E}^h_p(u)=\mathcal{\mu}^p(M, \partial{M}, h)$, then Corollary \ref{rigidityminimizer} implies
$$
u\in Fix(T_h).
$$
Using \eqref{deffix}, we get
$$
u=T_h(u)
$$
So, \eqref{sameval} gives
$$
\mathcal{I}^h_p(u)=\mathcal{E}^h_p(u).
$$
Therefore, using Proposition \ref{sameinf}, we get $$\mathcal{I}^h_p(u)=\mathcal{\mu}^p_{oc}(M, \partial{M}, h).$$
Similarly, if $\mathcal{I}^h_p(u)=\mathcal{\mu}^p_{oc}(M, \partial{M}, h)$, then as before, Corollary \ref{rigidityminimizeri} implies
$$
u=T_h(u).
$$
So, \eqref{sameval} gives
$$
\mathcal{E}^h_p(u)=\mathcal{I}^h_p(u).
$$
Therefore, using Proposition \ref{sameinf}, we obtain $$\mathcal{E}^h_p(u)=\mathcal{\mu}^p(M, \partial{M}, h),$$ and this ends the proof.
\end{pf}

\noindent
\begin{pfn}{ of Theorem \ref{positivevarmass}}\\
By Proposition \ref{sameminimizer} we have
\begin{equation}\label{imp00}
\mathcal{E}^g(u)=\mathcal{\mu}(M, \partial{M}, [g]).
\end{equation}
Thus as a critical point of $\mathcal{E}^g$, the strong maximum principle and the smooth regularity of positive weak solutions of the Cherrier-Escobar equation \eqref{eq:sequation} imply
\begin{equation}\label{imp01}
u\in C^{\infty}_+(\overline{M}, g)\quad \text{and} \quad
\begin{cases}
R_{g_{u}} = 0, \\
H_{g_{u}} = \text{const}>0.
\end{cases}
\end{equation}
Corollary  \ref{rigidityminimizer}, \eqref{deffix}, and \eqref{imp00} imply
\begin{equation}\label{imp02}
u\in Fix(T_g)\Longleftrightarrow u=T_g(u).
\end{equation}
Hence, point 1) follows from \eqref{imp01} and \eqref{imp02}.\\
Since $\mathcal{I}^g(u_{min}) = \mathcal{\mu}_{oc}(M, \partial{M}, [g])$ then $\exists \; \lambda>0$ such that by Proposition \ref{sameminimizer}, Corollary \ref{scaleinv}, and \eqref{phjph}, $u^{min} := \lambda u_{min}$ verifies
$$
\begin{cases}
R_{g_{u^{\min}}} = 0, \\
H_{g_{u^{\min}}} = \mu(M, \partial M, [g]),
\end{cases}
$$
and
\begin{equation}\label{imp1}
\mathcal{E}^g(u^{\min}) = \mu(M, \partial{M}, [g]).
\end{equation}
Thus Proposition \ref{sameminimizer} and \eqref{imp1} imply
\begin{equation}\label{imp3}
\mathcal{I}^g(u^{min})=\mathcal{\mu}_{oc}(M, \partial{M}, [g]).
\end{equation}
Hence, point 2) follows from point 1).\\
\newpage
\noindent
For 3), we argue as follows. Since
\begin{equation}
\mu_{oc}(M, \partial{M}, [g]) = \mu(M, \partial{M}, [g]) \tag{$a'$}
\end{equation}
see Proposition \ref{sameinf}, then $\mathcal{E}^g$ and $\mathcal{I}^g$ have the same minimizing sequences. Now, by Ekeland's variational principle \cite{ivek1}, \;$\mathcal{E}^g$ has a minimizing sequence which is a Palais-Smale (PS) sequence. Hence there exists a sequence $u_k \in H^1_+(\overline{M}, g)$ such that:\\\\
\textcircled{A}$\mathcal{I}^g(u_k) \to \mathcal{\mu}_{oc}(M, \partial{M}, [g])$,\\\\
\textcircled{B} $\mathcal{E}^g(u_k) \to \mathcal{\mu}(M, \partial{M}, [g])$,\\\\
\textcircled{C} $\nabla \mathcal{E}^g(u_k) \to 0 \text{ in }  H^1(\overline{M}, g)$.\\\\
By the work of Almaraz \cite{alm3}, \textcircled{B} and $\textcircled{C}$ implies that up to a subsequence and scaling, $u_k \rightharpoonup u$ in $H^1(\overline{M}, g)$,\;$u\geq0$ in $M$,\;$u\geq0$ on $\partial{M}$ and $u_k$ verifies the following alternatives. Either
$$(s1) \;\;\;  u_k \to u \text{ in } H^1_+(\overline{M}, g) \text{ or }$$
$$(s2) \;\;\;  \lim_{k \to \infty} \mathcal{E}^g(u_k) \geq \mathcal{\mu}(\mathbb{B}^n, \mathbb{S}^{n-1}, [g_{\overline{\mathbb{B}}^n}])$$
If $(s2)$ is true, then $(2)$ implies that 
\begin{equation}
\mu(M, \partial{M}, [g]) \geq \mu(\mathbb{B}^n, \mathbb{S}^{n-1}, [g_{\overline{\mathbb{B}}^n}])  \tag{$a''$}   
\end{equation}
On the other hand, since $(a')$ is true for any compact Riemannian manifold $(\overline{M}, g),$ then we have \begin{equation}
\mu(\mathbb{B}^n, \mathbb{S}^{n-1}, [g_{\overline{\mathbb{B}}^n}]) = \mu_{\text{oc}}(\mathbb{B}^n, \mathbb{S}^{n-1}, [g_{\overline{\mathbb{B}}^n}]) \tag{$a'''$}
\end{equation} 
Using $(a')$, $(a'')$ and $(a''')$ then we get $$\mu_{oc}(M, \partial{M}, [g]) \geq \mu_{oc}(\mathbb{B}^n, \mathbb{S}^{n-1}, [g_{\overline{\mathbb{B}}^n}],)$$ hence a contradiction to $(a)$. Hence $\exists \; u_{min} \in H^1_+(\overline{M}, g)$ such that up to a subsequence and scaling $$u_k \to u_{min} \in H^1(\overline{M}, g) \implies \lVert u_k \rVert^2_g \to \lVert u_{min} \rVert^2_g, \text{ and } \text{tr}(u_k) \to \text{tr}(u_{min}) \in L^{2^\#}(\partial{M}, \hat{g}).$$
\noindent
So we get that 
$$\mathcal{I}^g(u_k) \to \mathcal{I}^g(u_{min}).$$
Hence \textcircled{A}$ \implies \mathcal{I}^g(u_{min}) = \mathcal{\mu}_{oc}(M, \partial{M}, [g])$ thereby ending the proof.
\end{pfn}

\section{Obstacle problem and Sobolev trace type inequality}\label{casesphere}
In this section, we discuss some Sobolev trace type inequalities related to the boundary obstacle problem for the couple conformal Laplacian and conformal Robin operator. In particular, we specialize to the case of the \;$n$-dimensional standard ball \;$(\overline{\mathbb{B}}^n, g_{\overline{\mathbb{B}}^n})$.
\vspace{7pt}

\noindent
Lemma \ref{sobolev1} (or Proposition \ref{sameinf}) implies the following obstacle Sobolev trace type inequality.
\begin{lem}\label{sobolev1op}
Assuming $\mathcal{\mu}(M, \partial{M},  [g])>0$\; and \;$h\in[g]$, then for \;$u\in H^1_+(\overline{M}, h)$
$$
||T_h(u)||_{L^{2^\#}(\partial{M}, \hat{h})}\leq \frac{1}{\sqrt{\mathcal{\mu}(M, \partial{M}, [g])}}||u||_h
$$
\end{lem}
\begin{pf}
Lemma \ref{sobolev1} implies
$$
||T_h(u)||_{L^{2^\#}(\partial{M}, \hat{h})}\leq \frac{1}{\sqrt{\mathcal{\mu}(M,\partial{M}, [g])}}||T_h(u)||_h
$$
By the minimality property of $T_h(u)$, we obtain
$$
||T_h(u)||_h\leq ||u||_h.
$$
Thus combining the two inequalities, we get
$$
||T_h(u)||_{L^{2^\#}(\partial{M}, \hat{h})}\leq \frac{1}{\sqrt{\mathcal{\mu}(M, \partial{M}, [g])}}||u||_h
$$
as desired. Also, this result can be obtained using Proposition \ref{sameinf}. Indeed, using the definition of \;$\mathcal{I}^h$\; and \;$\mathcal{\mu}_{oc}(M, \partial{M}, [g])$ (see \eqref{eq:defi1} and Definition \ref{defyp}), and Proposition \ref{sameinf}, we have
$$
||T_h(u)||_{L^{2^\#}(\partial{M}, \hat{h})}\leq \frac{1}{\sqrt{\mathcal{\mu}_{oc}(M, \partial{M}, [g])}}||u||_h=\frac{1}{\sqrt{\mathcal{\mu}(M, \partial{M}, [g])}}||u||_h.
$$
\end{pf}
\vspace{7pt}

\noindent
When $(\overline{M}, g)=(\mathbb{B}^n, g_{\mathbb{B}^n})$, we have the following boundary obstacle analogue of Lemma \ref{sobolev2}.
\begin{pro}\label{sobolev2op}
Assuming $(\overline{M}, g)=(\overline{\mathbb{B}}^n, g_{\overline{\mathbb{B}}^n})$ and $h=g_w\in [g]$, then for $u\in H^1_+(\overline{M}, h)$,
\begin{equation}\label{inequality}
\|T_h(u)\|_{L^{2^\#}(\partial{M}, \hat{h})}\leq \frac{1}{\sqrt{\mathcal{\mu}(\mathbb{B}^n, \mathbb{S}^{n-1}, [g_{\overline{\mathbb{B}}^n}])}}\|u\|_h,
\end{equation}
with equality holdsing if and only if $u\in C^{\infty}_+(\overline{M}, h)$ and
\[
\begin{cases}
R_{g_{wu}} = 0, \\
H_{g_{wu}} = \text{const on}\;\;\;\partial{M}.
\end{cases}
\]
\end{pro}

\begin{pf}
The inequality \eqref{inequality} clearly follows from Lemma \ref{sobolev1op}, while the rigidity part follows from Proposition \ref{sameinf}, Proposition \ref{sameminimizer} and Lemma \ref{sobolev2}. Indeed, $u$ attains equality in \eqref{inequality} is equivalent to 
\[
\|T_h(u)\|_{L^{2^\#}(\partial{M}, \hat{h})}=\frac{1}{\sqrt{\mathcal{\mu}(\mathbb{B}^n, \mathbb{S}^{n-1},[g_{\overline{\mathbb{B}}^n}])}}\|u\|_h. 
\]
So, Proposition \ref{sameinf} implies $u$ attains equality in \eqref{inequality} is equivalent to 
\[
\|T_h(u)\|_{L^{2^\#}(\partial{M}, \hat{h})}=\frac{1}{\sqrt{\mathcal{\mu}_{oc}(\mathbb{B}^n, \mathbb{S}^{n-1}, [g_{\overline{\mathbb{B}}^n}])}}\|u\|_h.
\]
And using the definition of $\mathcal{I}^h$ (see \eqref{eq:defi1}), we have $u$ attains equality in \eqref{inequality} is equivalent to 
\[
\mathcal{I}^h(u)=\mathcal{\mu}_{oc}(\mathbb{B}^n, \mathbb{S}^{n-1}, [g_{\overline{\mathbb{B}}^n}]).
\] 
Furthermore, using Proposition \ref{sameminimizer}, we infer that $u$ attains equality in \eqref{inequality} is equivalent to 
\[
\mathcal{E}^h(u)=\mathcal{\mu}(\mathbb{B}^n, \mathbb{S}^{n-1}, [g_{\overline{\mathbb{B}}^n}]).
\]
Therefore, using the definition of $\mathcal{E}^h$ (see \eqref{eq:defj1}), we have $u$ attains equality in \eqref{inequality} is equivalent to $u$ attains equality in \eqref{inequality0}. Hence the rigidity part in Lemma \ref{sobolev2} implies $u$ attains equality in \eqref{inequality} is equivalent to $u\in C^{\infty}_+(\overline{M}, h)$ and
\[
\begin{cases}
R_{g_{wu}} = 0, \\
H_{g_{wu}} = \text{const on}\;\;\;\partial{M}.
\end{cases}
\]
\end{pf}

\vspace{7pt}

\noindent
\begin{pfn}{ of Theorem \ref{sphere}}\\
Theorem \ref{sphere}  follows directly from Proposition \ref{sameinf}, and the rigidity part in Proposition \ref{sobolev2op}. Indeed, Proposition \ref{sameinf} implies
$$\mathcal{I}^g(u)=\mathcal{\mu}_{oc}(M, \partial{M}, [g]) \Longleftrightarrow \mathcal{I}^g(u)=\mathcal{\mu}(M, \partial{M},  [g]).$$ 
\noindent
Therefore, by the definition of $\mathcal{I}^g$, we have $\mathcal{I}^g(u)=\mathcal{\mu}_{oc}(M, \partial{M}, [g])$ is equivalent to  $u$ attains equality in \eqref{inequality} with $g=h$. Hence, using the rigidity part of Proposition \ref{sobolev2op}, we have $\mathcal{I}^g(u)=\mathcal{\mu}_{oc}(M, \partial{M}, [g])$ is equivalent to $u\in C^{\infty}_+(\overline{M}, g)$ and
\[
\begin{cases}
R_{g_{u}} = 0, \\
H_{g_{u}} = \text{const on}\;\;\;\partial{M}.
\end{cases}
\]
\end{pfn}

\end{document}